# ELEMENTARY BOUNDS ON POINCARÉ AND LOG-SOBOLEV CONSTANTS FOR DECOMPOSABLE MARKOV CHAINS

By Mark Jerrum[1], Jung-Bae Son[2], Prasad Tetali[3] and Eric Vigoda[4]

*University of Edinburgh, University of Edinburgh, Georgia Institute of Technology and University of Chicago*

We consider finite-state Markov chains that can be naturally decomposed into smaller "projection" and "restriction" chains. Possibly this decomposition will be inductive, in that the restriction chains will be smaller copies of the initial chain. We provide expressions for Poincaré (resp. log-Sobolev) constants of the initial Markov chain in terms of Poincaré (resp. log-Sobolev) constants of the projection and restriction chains, together with further a parameter. In the case of the Poincaré constant, our bound is always at least as good as existing ones and, depending on the value of the extra parameter, may be much better. There appears to be no previously published decomposition result for the log-Sobolev constant. Our proofs are elementary and self-contained.

**1. The setting.** In a number of applications, one is interested in finding tight, nonasymptotic upper bounds on the mixing time, that is, rate of convergence to stationarity, of finite-state Markov chains. One important example arises in the analysis of Markov chain Monte Carlo algorithms. These are algorithms for sampling and counting combinatorial structures that are founded on Markov chain simulation. The efficiency of these algorithms depends crucially on the rate of convergence to stationarity of the Markov chain being simulated.

Received March 2003; revised September 2003.
[1]Supported in part by an EPSRC grant and the IST Programme of the EU under contract IST-1999-14036 (RAND-APX).
[2]Supported by an award from DAAD.
[3]Supported in part by NSF Grant DMS-01-00298.
[4]Supported in part by an NSF CAREER grant.

In proving rapid convergence to stationarity, Poincaré and latterly log-Sobolev inequalities have proved to be powerful tools. The larger the constants in these inequalities, the faster the convergence to stationarity. (These and other informal remarks are made rigorous in the following section.) Here we consider finite-state Markov chains whose description suggests a natural state-space partition. This partition naturally induces a number of restriction chains, in which transitions are restricted to occur within blocks of the partition, and a projection chain, whose states are the blocks themselves. The hope is that by computing Poincaré or log-Sobolev constants for the restriction and projection chains we can obtain a Poincaré or log-Sobolev constant for the original chain. Various authors, including Madras, Martin and Randall [[13]–[15]], have investigated this approach.

Sometimes it may be possible to apply the decomposition step inductively, as was done by Cooper, Dyer, Frieze and Rue [4] in the context of spin models on "narrow grids" and by Jerrum and Son [9] for the "bases-exchange walk" for balanced matroids. In these applications, it is particularly important that our arguments give as little as possible away at each decomposition step.

Clearly there is a need for general decomposition theorems that relate, say, the Poincaré constant $\lambda$ of the original chain as tightly as possible to those of the restriction and projection chains. The existing decomposition theorems of this sort seem all to rest ultimately on an unpublished result of Caracciolo, Pelissetto and Sokal [1]. (Note, however, that a statement of their result and a version of their proof were published as an appendix to an article by Madras and Randall [13].)

Our first goal, then, is to provide an elementary, self-contained and accessible account of the basic decomposition result. However, in developing the result from first principles we find we can derive a statement that is considerably sharper than the current ones in many situations. For example, in the context of inductively defined Markov chains, existing decomposition results cannot yield inverse polynomial bounds on $\lambda$ (and hence polynomial bounds on mixing time), even when the depth of the induction is logarithmic. In contrast, we are able to give inverse polynomial bounds on $\lambda$ for inductively defined Markov chains, and are even able to recover the result of Jerrum and Son [9] on the bases-exchange walk for balanced matroids, where the depth of the induction is linear in some natural measure of the input size.

It transpires that the proof of the decomposition result for the Poincaré constant carries over straightforwardly to the log-Sobolev constant. In many situations the optimal log-Sobolev constant seems to be within a small constant factor of the optimal Poincaré constant (spectral gap); the advantage of the log-Sobolev constant in these situations is that it translates to a tighter bound on mixing time (construed as the time to convergence to near stationarity in $\ell_1$ norm). To the best of our knowledge, this is the first general decomposition result for log-Sobolev inequalities, although it should be mentioned that Cesi [2] gave an argument that applies when the state space is a Cartesian product.

We have stated our results for finite-state Markov chains, since that seems to be the natural setting given the potential applications to Markov chain Monte Carlo algorithms. However, everything extends (with no notational change) to countably infinite state spaces and (with appropriate notational changes and possible regularity conditions) to uncountable state spaces.

**2. Poincaré constant via decomposition.** Consider an ergodic Markov chain on finite state space $\Omega$ with transition probabilities $P:\Omega \times \Omega \to [0,1]$ and stationary distribution $\pi:\Omega \to [0,1]$. We assume that the Markov chain is *time-reversible*, that is to say, it satisfies the *detailed balance* condition

$$\pi(x)P(x,y) = \pi(y)P(y,x) \qquad \text{for all } x,y \in \Omega.$$

Let $\Omega = \Omega_0 \cup \cdots \cup \Omega_{m-1}$ be a decomposition of the state space into $m$ disjoint sets. As usual, we use $[m] := \{0, \ldots, m-1\}$ to denote the first $m$ natural numbers.

Following [15], we define $\bar{\pi}:[m] \to [0,1]$ by

$$\bar{\pi}(i) = \sum_{x \in \Omega_i} \pi(x)$$

and define $\overline{P}:[m] \times [m] \to [0,1]$ by

$$\overline{P}(i,j) := \bar{\pi}(i)^{-1} \sum_{\substack{x \in \Omega_i \\ y \in \Omega_j}} \pi(x)P(x,y).$$

The Markov chain on state space $[m]$ and with transition probabilities $\overline{P}$ is the *projection* chain induced by the partition $\{\Omega_i\}$. Since the original Markov chain is time-reversible, so is the projection chain. It is easy to check, using this observation, that the projection chain has $\bar{\pi}$ as a stationary distribution.

For each $i \in [m]$ the *restriction* Markov chain on $\Omega_i$ has transition probabilities $P_i:\Omega_i \times \Omega_i \to [0,1]$ defined by

$$P_i(x,y) = \begin{cases} P(x,y), & \text{if } x \neq y, \\ 1 - \sum_{z \in \Omega_i \setminus \{x\}} P(x,z), & \text{if } x = y. \end{cases}$$

Again, the restriction chain inherits time-reversibility from the original chain, and so it has $\pi_i:\Omega_i \to [0,1]$ as a stationary distribution, where $\pi_i(x) = \pi(x)/\bar{\pi}(i)$. In applications, we require the projection chain and all the restriction chains to be irreducible; in which case the various stationary distributions $\bar{\pi}$ and $\pi_0, \ldots, \pi_{m-1}$ are unique.

Let $f:\Omega \to \mathbb{R}$ be an arbitrary test function. The expectation and variance of $f$ with respect to $\pi$ are of course given by

$$\mathrm{E}_\pi f := \sum_{x \in \Omega} \pi(x) f(x)$$

and

$$\mathrm{Var}_\pi f := \sum_{x \in \Omega} \pi(x)(f(x) - \mathrm{E}_\pi f)^2,$$

respectively. The *Dirichlet form* associated with $f$ and $P$ is defined as

$$\mathcal{E}_\pi(f,f) := \tfrac{1}{2} \sum_{x,y \in \Omega} \pi(x) P(x,y) (f(x) - f(y))^2.$$

Consider now a *Poincaré inequality* of the form

$$\mathcal{E}_\pi(f,f) \geq \lambda \, \mathrm{Var}_\pi f \tag{1}$$

that holds uniformly over all $f:\Omega \to \mathbb{R}$, with $\lambda > 0$ being the corresponding *Poincaré constant*. It is well known that a lower bound on $\lambda$ translates directly to an upper bound on mixing time of a Markov chain. To avoid technical problems associated with nearly periodic Markov chains, assume that loop probabilities are uniformly bounded away from 0. [Alternatively, interpret $P(\cdot,\cdot)$ as the transition rates of a continuous-time Markov chain on $\Omega$.] Denote by $P^t(x,\cdot)$ the $t$-step distribution of the chain, given that $x \in \Omega$ is the initial state. Then there is a function $t:\Omega \times (0,1] \to \mathbb{N}$ with

$$t(x,\varepsilon) = O\!\left(\frac{1}{\lambda}\!\left(\ln \frac{1}{\pi(x)} + \ln \frac{1}{\varepsilon}\right)\right) \tag{2}$$

such that $\|P^{t(x,\varepsilon)}(x,\cdot) - \pi\|_{\mathrm{TV}} \leq \varepsilon$, where $\|\cdot\|_{\mathrm{TV}}$ is total variation norm (i.e., half the $\ell_1$ norm). For a proof of this claim that is valid for general (i.e., not necessarily time-reversible) Markov chains, refer to [8], Section 5.2, interpreting $\varrho$ there as the reciprocal of our $\lambda$.

Observe that in our notation for expectation, variance and so forth, we make explicit the probability distribution $\pi$ as a subscript, because this varies throughout our proofs. For example, we may write Poincaré inequalities for the projection and restriction chains as $\mathcal{E}_{\bar{\pi}}(\bar{f},\bar{f}) \geq \bar{\lambda} \, \mathrm{Var}_{\bar{\pi}} \bar{f}$ and $\mathcal{E}_{\pi_i}(f,f) \geq \lambda_i \, \mathrm{Var}_{\pi_i} f$, respectively. Naturally, $\bar{\pi}$ (resp. $\pi_i$) is to be considered as a probability distribution on $[m]$ (resp. $\Omega_i$), and $\bar{f}$ as a function on $[m]$.

Suppose the projection chain and the various restriction chains satisfy Poincaré inequalities with constants $\bar{\lambda}$, and $\lambda_0, \ldots, \lambda_{m-1}$, respectively. Define $\lambda_{\min} = \min_i \lambda_i$. We are interested in obtaining a Poincaré inequality for the original Markov chain, with Poincaré constant $\lambda = \lambda(\bar{\lambda}, \lambda_{\min}, \gamma)$, where $\gamma$ is a further parameter

(3) $$\gamma := \max_{i \in [m]} \max_{x \in \Omega_i} \sum_{y \in \Omega \setminus \Omega_i} P(x, y).$$

Of course, we would like $\lambda$ to be as large as possible. Informally, $\gamma$ is the probability of escape in one step from the current block of the partition, maximized over all states. Given this interpretation, it is clear that $\gamma$ never exceeds 1, and may be much smaller in many applications. It is in these applications that we improve on existing decomposition bounds.

THEOREM 1. *Consider a finite-state time-reversible Markov chain decomposed into a projection chain and $m$ restriction chains as above. Suppose the projection chain satisfies a Poincaré inequality with constant $\bar{\lambda}$, and the restriction chains satisfy inequalities with uniform constant $\lambda_{\min}$. Let $\gamma$ be defined as in* (3). *Then the original Markov chain satisfies a Poincaré inequality with constant*

$$\lambda := \min\left\{\frac{\bar{\lambda}}{3}, \frac{\bar{\lambda}\lambda_{\min}}{3\gamma + \bar{\lambda}}\right\}.$$

PROOF. Consider an arbitrary test function $f : \Omega \to \mathbb{R}$. Our starting point is the decomposition of $\operatorname{Var}_\pi f$ with respect to the partition $\Omega = \Omega_0 \cup \cdots \cup \Omega_{m-1}$,

(4) $$\operatorname{Var}_\pi f = \sum_{i \in [m]} \bar{\pi}(i) \operatorname{Var}_{\pi_i} f + \sum_{i \in [m]} \bar{\pi}(i)(\operatorname{E}_{\pi_i} f - \operatorname{E}_\pi f)^2,$$

and a similar decomposition of the Dirichlet form,

(5) $$\mathcal{E}_\pi(f, f) = \sum_{i \in [m]} \bar{\pi}(i) \mathcal{E}_{\pi_i}(f, f) + \tfrac{1}{2} \sum_{\substack{i,j \in [m] \\ i \neq j}} \mathcal{C}_{ij},$$

where

$$\mathcal{C}_{ij} := \sum_{\substack{x \in \Omega_i \\ y \in \Omega_j}} \pi(x) P(x, y)(f(x) - f(y))^2.$$

Identity (5) is almost content-free and comes from partitioning the terms in the definition of $\mathcal{E}_\pi(f, f)$ according to whether $i$ and $j$ are in the same or in different blocks of the partition. Identity (4) has a little more substance, but is nevertheless standard and can be obtained through simple algebraic manipulation. It states informally that the variance of $f$ may be obtained by summing the variances *within* blocks of the partition and the variance *between* blocks.

In summations and so forth, variables $i$ and $j$ always range over $[m]$, so we are not explicit about their range in what follows. For all $i,j$ with $i \neq j$ and $\overline{P}(i,j) > 0$, define $\hat{\pi}_i^j : \Omega_i \to [0,1]$ by

$$\hat{\pi}_i^j(x) := \frac{\pi_i(x) \sum_{y \in \Omega_j} P(x,y)}{\overline{P}(i,j)}.$$

Note that $\hat{\pi}_i^j$ is a probability distribution on $\Omega_i$.

The first term on the right-hand side of (4) we simply bound as

$$\sum_i \bar{\pi}(i) \operatorname{Var}_{\pi_i} f \leq \sum_i \frac{1}{\lambda_i} \bar{\pi}(i) \mathcal{E}_{\pi_i}(f,f)$$

(6)

$$\leq \frac{1}{\lambda_{\min}} \sum_i \bar{\pi}(i) \mathcal{E}_{\pi_i}(f,f).$$

The second term we transform, starting with an application of the Poincaré inequality for the projection chain,

$$\sum_i \bar{\pi}(i)(\mathrm{E}_{\pi_i} f - \mathrm{E}_\pi f)^2$$

$$\leq \frac{1}{2\bar{\lambda}} \sum_{i \neq j} \bar{\pi}(i) \overline{P}(i,j)(\mathrm{E}_{\pi_i} f - \mathrm{E}_{\pi_j} f)^2$$

$$\leq \frac{3}{2\bar{\lambda}} \sum_{i \neq j} \bar{\pi}(i) \overline{P}(i,j)$$

$$\times [(\mathrm{E}_{\pi_i} f - \mathrm{E}_{\hat{\pi}_i^j} f)^2 + (\mathrm{E}_{\hat{\pi}_i^j} f - \mathrm{E}_{\hat{\pi}_j^i} f)^2 + (\mathrm{E}_{\hat{\pi}_j^i} f - \mathrm{E}_{\pi_j} f)^2]$$

(7) $\qquad = \frac{3}{2\bar{\lambda}}[\Sigma_1 + \Sigma_2 + \Sigma_3],$

where

$$\Sigma_1 := \sum_{i \neq j} \bar{\pi}(i) \overline{P}(i,j)(\mathrm{E}_{\pi_i} f - \mathrm{E}_{\hat{\pi}_i^j} f)^2, \ldots.$$

We proceed to bound $\Sigma_1$, $\Sigma_2$ and $\Sigma_3$ separately, noting that $\Sigma_1$ and $\Sigma_3$ are equal by time-reversibility. For the second of these we have

(8) $\qquad \Sigma_2 = \sum_{i \neq j} \bar{\pi}(i) \overline{P}(i,j) \left[ \sum_{\substack{x \in \Omega_i \\ y \in \Omega_j}} \frac{\pi(x) P(x,y)}{\bar{\pi}(i) \overline{P}(i,j)}(f(x) - f(y)) \right]^2$

(9) $\qquad \leq \sum_{i \neq j} \bar{\pi}(i) \overline{P}(i,j) \sum_{\substack{x \in \Omega_i \\ y \in \Omega_j}} \frac{\pi(x) P(x,y)}{\bar{\pi}(i) \overline{P}(i,j)}(f(x) - f(y))^2$

$$= \sum_{i \neq j} \sum_{\substack{x \in \Omega_i \\ y \in \Omega_j}} \pi(x) P(x,y)(f(x) - f(y))^2$$

$$(10) \qquad = \sum_{i \neq j} \mathcal{C}_{ij},$$

where (8) uses the fact that $\pi(x)P(x,y)/\bar{\pi}(i)\overline{P}(i,j)$ is a joint distribution on $\Omega_i \times \Omega_j$ whose marginals are $\hat{\pi}_i^j$ and $\hat{\pi}_j^i$, and (9) is seen to be Cauchy–Schwarz, once we have noted that

$$\sum_{\substack{x \in \Omega_i \\ y \in \Omega_j}} \frac{\pi(x) P(x,y)}{\bar{\pi}(i) \overline{P}(i,j)} = 1,$$

by definition.

Now for $\Sigma_1$. Using standard facts about variance,

$$\operatorname{Var}_{\hat{\pi}_i^j} f = \operatorname{Var}_{\hat{\pi}_i^j}(f - \mathrm{E}_{\pi_i} f)$$

$$(11) \qquad = \sum_{x \in \Omega_i} \hat{\pi}_i^j(x)(f(x) - \mathrm{E}_{\pi_i} f)^2 - (\mathrm{E}_{\hat{\pi}_i^j} f - \mathrm{E}_{\pi_i} f)^2,$$

so that certainly

$$(12) \qquad (\mathrm{E}_{\hat{\pi}_i^j} f - \mathrm{E}_{\pi_i} f)^2 \leq \sum_{x \in \Omega_i} \hat{\pi}_i^j(x)(f(x) - \mathrm{E}_{\pi_i} f)^2.$$

Thus we have the bound

$$\Sigma_1 \leq \sum_{i \neq j} \bar{\pi}(i)\overline{P}(i,j) \sum_{x \in \Omega_i} \hat{\pi}_i^j(x)(f(x) - \mathrm{E}_{\pi_i} f)^2$$

$$= \sum_i \bar{\pi}(i) \sum_{x \in \Omega_i} \pi_i(x)(f(x) - \mathrm{E}_{\pi_i} f)^2 \sum_{j : j \neq i} \frac{\hat{\pi}_i^j(x) \overline{P}(i,j)}{\pi_i(x)}$$

$$(13) \qquad = \sum_i \bar{\pi}(i) \sum_{x \in \Omega_i} \pi_i(x)(f(x) - \mathrm{E}_{\pi_i} f)^2 \sum_{j : j \neq i} P(x, \Omega_j)$$

$$(14) \qquad \leq \gamma \sum_i \bar{\pi}(i) \operatorname{Var}_{\pi_i} f$$

$$(15) \qquad \leq \frac{\gamma}{\lambda_{\min}} \sum_i \bar{\pi}(i) \mathcal{E}_{\pi_i}(f,f),$$

where (13) applies the definition of $\hat{\pi}_i^j$, (14) applies the definition of $\gamma$ and (15) applies the Poincaré inequalities for the restriction chains.

Substituting (10) and (15) in (7), and recalling that $\Sigma_1 = \Sigma_3$, we have

$$(16) \quad \sum_i \bar{\pi}(i)(\mathrm{E}_{\pi_i} f - \mathrm{E}_\pi f)^2 \leq \frac{3}{2\bar{\lambda}} \sum_{i \neq j} \mathcal{C}_{ij} + \frac{3\gamma}{\bar{\lambda}\lambda_{\min}} \sum_i \bar{\pi}(i) \mathcal{E}_{\pi_i}(f,f).$$

Then substituting (6) and (16) into (4) yields

$$\text{(17)} \quad \text{Var}_\pi f \leq \frac{3}{2\bar{\lambda}} \sum_{i \neq j} \mathcal{C}_{ij} + \frac{3\gamma + \bar{\lambda}}{\bar{\lambda}\lambda_{\min}} \sum_i \bar{\pi}(i)\,\mathcal{E}_{\pi_i}(f,f).$$

Finally, comparing (17) with (5), we see that

$$\mathcal{E}_\pi(f,f) \geq \lambda \, \text{Var}_\pi f,$$

where $\lambda$ is as in the statement of the theorem. $\square$

The first thing to note is that $\gamma \leq 1$, so that always $\lambda = \Omega(\bar{\lambda}\lambda_{\min})$, matching existing decomposition results (e.g., [1]). It may be the case that $\gamma$ is indeed a constant (e.g., decompose a random walk on $[n]$ into $k$ random walks on $[n/k]$, where we assume for convenience that $k$ divides $n$). In this case, $\gamma = \frac{1}{2}$ and we get no improvement over existing bounds.

At the other extreme, there are situations, for example, spin systems on fragments of the Bethe lattice or narrow grids, where $\gamma$ and $\bar{\lambda}$ are of the same order of magnitude. Applying Theorem 1 inductively then yields bounds on the spectral gap that are inverse polynomial in the problem size $n$, provided the depth of recursion is $O(\log n)$. Section 4.3 treats such an example.

This seems about the best that can be achieved using a parameter as "global" as $\gamma$. To go further, we need a much stricter pointwise constraint on the distributions $\hat{\pi}_i^j$. For example, if we know that

$$\text{(18)} \quad (1-\eta)\pi_i \leq \hat{\pi}_i^j \leq (1+\eta)\pi_i,$$

pointwise, whenever $\hat{\pi}_i^j$ is defined [i.e., whenever $\overline{P}(i,j) > 0$], then

$$\text{(19)} \quad (\text{E}_{\pi_i} f - \text{E}_{\hat{\pi}_i^j} f)^2 \leq (1+\eta)\text{Var}_{\pi_i} f - \text{Var}_{\hat{\pi}_i^j} f$$

$$\leq (1+\eta)\text{Var}_{\pi_i} f - (1-\eta) \sum_{x \in \Omega_i} \pi_i(x)(f(x) - \text{E}_{\hat{\pi}_i^j} f)^2$$

$$\text{(20)} \quad \leq (1+\eta)\text{Var}_{\pi_i} f - (1-\eta)\text{Var}_{\pi_i} f$$

$$\text{(21)} \quad = 2\eta \, \text{Var}_{\pi_i} f,$$

where (19) comes from (11), and (20) from the fact that $\sum_{x \in \Omega_i} \pi_i(x)(f(x) - \mu)^2$ is minimized at $\mu = \text{E}_{\pi_i} f$. Introducing the modified parameter

$$\text{(22)} \quad \hat{\gamma} := 2\eta \max_{i \in [m]} \sum_{j:\, j \neq i} \overline{P}(i,j)$$

we obtain:

COROLLARY 2. *Suppose that* (18) *is satisfied for some* $\eta > 0$ *and that* $\hat{\gamma}$ *is as defined as in* (22). *Then Theorem* 1 *holds with* $\hat{\gamma}$ *replacing* $\gamma$.

PROOF. Simply use (21) in place of (12) in the derivation of inequality (15). □

Note that $\hat{\gamma}$ may even be 0 (which happens if $\eta = 0$), as in the case of the $n$-dimensional Boolean cube. When that happens, $\Sigma_1 = \Sigma_3 = 0$ and we save a factor 3 in the argument, leading to:

COROLLARY 3. *If $\hat{\gamma} = 0$, then Theorem 1 holds with $\lambda := \min\{\bar{\lambda}, \lambda_{\min}\}$.*

For the Boolean cube, Corollary 3 immediately gives the exact bound on spectral gap. Even when $\hat{\gamma} > 0$ we may be able to compare the given Markov chain with one with reduced transition probabilities for which $\hat{\gamma} = 0$. For example, in the case of the bases-exchange walk on a balanced matroid, we may "thin down" the transition probabilities between $\Omega_i$ and $\Omega_j$ until they form a fractional matching (which is possible by a result in [7]). Thus we recover the known bound on the spectral gap for balanced matroids. All this is covered in detail in Section 4.5.

**3. Log-Sobolev constant via decomposition.** The program described above extends to the log-Sobolev constant with little work. Following Diaconis and Saloff-Coste [5] (and others), define the entropy-like quantity

$$\mathcal{L}_\pi(f) := \mathrm{E}_\pi[f^2(\ln f^2 - \ln(\mathrm{E}_\pi f^2))]. \tag{23}$$

Again, we indicate the probability distribution $\pi$ explicitly as a subscript, so we can talk about $\mathcal{L}_{\pi_i}(f)$ and so forth. A *log-Sobolev inequality* is an inequality of the form

$$\mathcal{E}_\pi(f, f) \geq \alpha \mathcal{L}_\pi(f)$$

that holds for all $f : \Omega \to \mathbb{R}$. The motivation for studying the *log-Sobolev constant* $\alpha$ is the analogue of (2),

$$t(x, \varepsilon) = O\left(\frac{1}{\alpha}\left(\ln\ln\frac{1}{\pi(x)} + \ln\frac{1}{\varepsilon}\right)\right), \tag{24}$$

which provides an estimate, this time in terms of $\alpha$, for the number of steps sufficient to achieve $\|P^{t(x,\varepsilon)}(x, \cdot) - \pi\|_{\mathrm{TV}} \leq \varepsilon$. [To avoid trivialities, assume $\pi(x) \leq e^{-1}$ in (24).] The estimate (24) of mixing time may be read off from [5], equation (3.3), assuming loop probabilities are bounded away from 0. (Diaconis and Saloff-Coste worked in continuous time, avoiding potential problems associated with nearly periodic Markov chains.) The key point to note is that $\ln(1/\pi(x))$ in (2) is replaced by $\ln\ln(1/\pi(x))$ in (24). This may not seem like a major improvement until we recall that $\pi(x)$ is typically exponentially small as a function of instance size.

Our aim, then, is to find $\alpha = \alpha(\bar{\alpha}, \alpha_{\min}, \gamma)$ that satisfies

$$\mathcal{E}_\pi(f, f) \geq \alpha \mathcal{L}_\pi(f).$$

Obviously we want $\alpha$ to be as large as possible. Our analogue of Theorem 1 is:

THEOREM 4. *The setting is exactly as in Theorem 1. Suppose the projection chain satisfies a log-Sobolev inequality with constant $\bar{\alpha}$ and that the restriction chains satisfy inequalities with uniform constant $\alpha_{\min}$. Then the original Markov chain satisfies a log-Sobolev inequality with constant*

$$\alpha := \min\left\{\frac{\bar{\alpha}}{3}, \frac{\bar{\alpha}\alpha_{\min}}{3\gamma + \bar{\alpha}}\right\}.$$

PROOF. Just as with variance, $\mathcal{L}_\pi(f)$ may be decomposed with respect to the partition $\Omega = \Omega_0 \cup \cdots \cup \Omega_{m-1}$, the analogue of identity (4):

$$(25) \quad \mathcal{L}_\pi(f) = \sum_i \bar{\pi}(i) \mathcal{L}_{\pi_i}(f) + \sum_i \bar{\pi}(i)(\mathrm{E}_{\pi_i} f^2)(\ln(\mathrm{E}_{\pi_i} f^2) - \ln(\mathrm{E}_\pi f^2)).$$

Since $\mathcal{L}_\pi(f)$ is a less familiar quantity than variance, we offer, in an addendum, a derivation not only of (25), but also of a number of other identities and inequalities used in this section. By analogy with (4), the first term expresses the entropy within blocks of the partition and the second term expresses the entropy between blocks. [Compare the second term in (25) with the right-hand side of (23), observing how $\sqrt{\mathrm{E}_{\pi_i} f^2}$ is now the appropriate "averaging of $f$" over $\Omega_i$, taking on the role of $\mathrm{E}_{\pi_i} f$ in the earlier calculation.] The decomposition of entropy expressed in (25) has been exploited by other authors (e.g., [12]).

We deal with the first term exactly as before:

$$(26) \quad \sum_i \bar{\pi}(i) \mathcal{L}_{\pi_i}(f) \leq \sum_i \frac{1}{\alpha_i} \bar{\pi}(i) \mathcal{E}_{\pi_i}(f, f)$$
$$\leq \frac{1}{\alpha_{\min}} \sum_i \bar{\pi}(i) \mathcal{E}_{\pi_i}(f, f).$$

The second term we transform, in an analogous manner to (7), starting with an application of the log-Sobolev inequality for the projection chain,

$$\sum_i \bar{\pi}(i)(\mathrm{E}_{\pi_i} f^2)(\ln(\mathrm{E}_{\pi_i} f^2) - \ln(\mathrm{E}_\pi f^2))$$
$$\leq \frac{1}{2\bar{\alpha}} \sum_{i \neq j} \bar{\pi}(i) \overline{P}(i,j)(\sqrt{\mathrm{E}_{\pi_i} f^2} - \sqrt{\mathrm{E}_{\pi_j} f^2})^2$$

$$\leq \frac{3}{2\bar{\alpha}} \sum_{i\neq j} \bar{\pi}(i)\overline{P}(i,j)[(\sqrt{\mathrm{E}_{\pi_i} f^2} - \sqrt{\mathrm{E}_{\hat{\pi}_i^j} f^2})^2$$

$$+ (\sqrt{\mathrm{E}_{\hat{\pi}_i^j} f^2} - \sqrt{\mathrm{E}_{\hat{\pi}_j^i} f^2})^2$$

(27)

$$+ (\sqrt{\mathrm{E}_{\hat{\pi}_j^i} f^2} - \sqrt{\mathrm{E}_{\pi_j} f^2})^2]$$

(28) $$= \frac{3}{2\bar{\alpha}}[\Sigma_1 + \Sigma_2 + \Sigma_3],$$

where

$$\Sigma_1 = \sum_{i\neq j} \bar{\pi}(i)\overline{P}(i,j)(\sqrt{\mathrm{E}_{\pi_i} f^2} - \sqrt{\mathrm{E}_{\hat{\pi}_i^j} f^2})^2, \ldots.$$

Tackling $\Sigma_2$ first, we have

(29) $$\Sigma_2 \leq \sum_{i\neq j} \bar{\pi}(i)\overline{P}(i,j) \sum_{\substack{x\in\Omega_i \\ y\in\Omega_j}} \frac{\pi(x)P(x,y)}{\bar{\pi}(i)\overline{P}(i,j)}(f(x)-f(y))^2$$

$$= \sum_{i\neq j} \sum_{\substack{x\in\Omega_i \\ y\in\Omega_j}} \pi(x)P(x,y)(f(x)-f(y))^2$$

(30) $$= \sum_{i\neq j} \mathcal{C}_{ij},$$

where (29) is Jensen's inequality applied to the convex function $g(\xi,\zeta) := (\sqrt{\xi} - \sqrt{\zeta})^2$ defined on $\xi,\zeta \geq 0$. (This is an artifice borrowed from [11]; see Section 3.1 for details and see [12] for related applications.)

For $\Sigma_1$ (which equals $\Sigma_3$ by time-reversibility), we have the bound

(31) $$\Sigma_1 \leq \sum_{i\neq j} \bar{\pi}(i)\overline{P}(i,j) \sum_{x,y\in\Omega_i} \pi_i(x)\hat{\pi}_i^j(y)(f(x)-f(y))^2$$

(32) $$= \sum_i \bar{\pi}(i) \sum_{x,y\in\Omega_i} \pi_i(x)\pi_i(y)(f(x)-f(y))^2 \sum_{j:j\neq i} P(y,\Omega_j)$$

(33) $$\leq 2\gamma \sum_i \bar{\pi}(i) \mathrm{Var}_{\pi_i} f$$

(34) $$\leq \frac{2\gamma}{\lambda_{\min}} \sum_i \bar{\pi}(i) \mathcal{E}_{\pi_i}(f,f)$$

(35) $$\leq \frac{\gamma}{\alpha_{\min}} \sum_i \bar{\pi}(i) \mathcal{E}_{\pi_i}(f,f),$$

where (31) recycles the Jensen artifice (see Section 3.1), (32) applies the definition of $\hat{\pi}_i^j$, (33) applies the definition of $\gamma$, (34) applies the Poincaré inequalities for the restriction chains and finally (35) applies the general inequality $\alpha \leq \lambda/2$ (see, e.g., [5]) that relates Poincaré and log-Sobolev constants. [Strictly speaking, we must interpret $\lambda_{\min}$ here as the minimum over the *optimal* Poincaré constants (i.e., spectral gaps) of the $m$ restriction Markov chains.]

Substituting (30) and (35) into (28), and recalling that $\Sigma_1 = \Sigma_3$, we have

$$\sum_i \bar{\pi}(i)(\mathrm{E}_{\pi_i} f^2)(\ln(\mathrm{E}_{\pi_i} f^2) - \ln(\mathrm{E}_\pi f^2))$$
(36)
$$\leq \frac{3}{2\bar{\alpha}} \sum_{i \neq j} \mathcal{C}_{ij} + \frac{3\gamma}{\bar{\alpha}\alpha_{\min}} \sum_i \bar{\pi}(i) \mathcal{E}_{\pi_i}(f,f).$$

Then substituting (26) and (36) into (25) yields

$$(37) \qquad \mathcal{L}_\pi(f) \leq \frac{3}{2\bar{\alpha}} \sum_{i \neq j} \mathcal{C}_{ij} + \frac{3\gamma + \bar{\alpha}}{\bar{\alpha}\alpha_{\min}} \sum_i \bar{\pi}(i) \mathcal{E}_{\pi_i}(f,f).$$

Finally, comparing (37) with (5), we see that $\mathcal{E}_\pi(f,f) \geq \alpha \mathcal{L}_\pi(f)$, where $\alpha$ is as in the statement of the theorem. $\square$

The remarks following the proof of Theorem 1 apply also to Theorem 4. In particular, if $(1-\eta)\pi_i \leq \hat{\pi}_i^j \leq (1+\eta)\pi_i$, then

$$(38) \qquad (\sqrt{\mathrm{E}_{\pi_i} f^2} - \sqrt{\mathrm{E}_{\hat{\pi}_i^j} f^2})^2 \leq 2\eta \,\mathrm{Var}_{\pi_i} f,$$

using Jensen's inequality again, but this time with an optimal coupling of the two random variables (r.v.'s) (see Section 3.1), yielding:

COROLLARY 5. *Suppose that* (18) *is satisfied for some* $\eta > 0$ *and that* $\hat{\gamma}$ *is as defined as in* (22). *Then Theorem* 4 *holds with* $\hat{\gamma}$ *replacing* $\gamma$.

Again, when $\hat{\gamma} = 0$ we save a factor 3.

COROLLARY 6. *If* $\hat{\gamma} = 0$, *then Theorem* 4 *holds with* $\alpha := \min\{\bar{\alpha}, \alpha_{\min}\}$.

3.1. *Addendum*: *proofs of an identity and some inequalities.* This addendum contains derivations of some of the possibly less obvious identities and inequalities used above.

PROOF OF IDENTITY (25). By appropriate scaling of the function $f$, it is enough to establish (25) when $\mathrm{E}_\pi f^2 = 1$. With this simplification,

$$\mathcal{L}_\pi(f) = \mathrm{E}_\pi[f^2 \ln f^2] = \sum_i \bar{\pi}(i) \mathrm{E}_{\pi_i}[f^2 \ln f^2]$$

and
$$\sum_i \bar{\pi}(i)\mathcal{L}_{\pi_i}(f) = \sum_i \bar{\pi}(i) \operatorname{E}_{\pi_i}[f^2(\ln f^2 - \ln(\operatorname{E}_{\pi_i} f^2))].$$

Subtracting the latter from the former, we obtain
$$\mathcal{L}_\pi(f) - \sum_i \bar{\pi}(i)\mathcal{L}_{\pi_i}(f) = \sum_i \bar{\pi}(i)(\operatorname{E}_{\pi_i} f^2)\ln(\operatorname{E}_{\pi_i} f^2),$$

as required. □

PROOF OF INEQUALITY (29). Let $X$ and $Y$ be r.v.'s taking values in $\Omega_i$ and $\Omega_j$, respectively, and with joint distribution given by
$$\Pr(X = x \wedge Y = y) = \frac{\pi(x)P(x,y)}{\bar{\pi}(i)\overline{P}(i,j)}.$$

Thus the marginal distribution of $X$ (resp. $Y$) is $\hat{\pi}_i^j$ (resp. $\hat{\pi}_j^i$). Since, by calculus, $g(\xi,\zeta) := (\sqrt{\xi} - \sqrt{\zeta})^2$ is convex in $\xi, \zeta \geq 0$, Jensen's inequality yields
$$\begin{aligned}(\sqrt{\operatorname{E}_{\hat{\pi}_i^j} f^2} - \sqrt{\operatorname{E}_{\hat{\pi}_j^i} f^2})^2 &= (\sqrt{\operatorname{E}[f(X)^2]} - \sqrt{\operatorname{E}[f(Y)^2]})^2 \\ &\leq \operatorname{E}[(f(X) - f(Y))^2] \\ &= \sum_{\substack{x \in \Omega_i \\ y \in \Omega_j}} \frac{\pi(x)P(x,y)}{\bar{\pi}(i)\overline{P}(i,j)}(f(x) - f(y))^2. \quad \square\end{aligned}$$

PROOF OF INEQUALITY (31). As above, but now with $X$ and $Y$ being *independent* r.v.'s with the appropriate distributions. □

PROOF OF INEQUALITY (38). Assume $\pi_i \neq \hat{\pi}_i^j$, otherwise there is nothing to demonstrate. Let $X$ and $Y$ be r.v.'s that take values in $\Omega_i$, with joint distribution satisfying the following conditions: (i) $X$ has distribution $\pi_i$, (ii) $Y$ has distribution $\hat{\pi}_i^j$ and (iii) $\Pr(X \neq Y) = \|\pi_i - \hat{\pi}_i^j\|_{\mathrm{TV}}$. It is well known that such an *optimal coupling* of two r.v.'s exists. Denote by
$$\psi(x,y) := \Pr(X = x \wedge Y = y)$$

the joint distribution of $X$ and $Y$. Define
$$\widehat{\psi}(x,y) := \begin{cases} 0, & \text{if } x = y, \\ \psi(x,y), & \text{otherwise.} \end{cases}$$

Partition $\Omega$ into two sets $\Omega = A \cup B$ such that $\pi_i(x) - \hat{\pi}_i^j(x) > 0$ for all $x \in A$, and $\hat{\pi}_i^j(y) - \pi_i(y) \geq 0$ for all $y \in B$. By assumption, $A$ and $B$ are non-empty. Optimality of the coupling of $X$ and $Y$ entails $\sum_y \widehat{\psi}(x,y) =$ $\max\{\pi_i(x) - \hat{\pi}_i^j(x), 0\}$ and $\sum_x \widehat{\psi}(x,y) = \max\{\hat{\pi}_i^j(y) - \pi_i(y), 0\}$. Thus, by Jensen's inequality,

$$
\begin{aligned}
(\sqrt{\mathrm{E}_{\pi_i} f^2} &- \sqrt{\mathrm{E}_{\hat{\pi}_i^j} f^2})^2 \\
&= (\sqrt{\mathrm{E}[f(X)^2]} - \sqrt{\mathrm{E}[f(Y)^2]})^2 \\
&\leq \mathrm{E}[(f(X) - f(Y))^2] \\
&= \sum_{x,y \in \Omega_i} \widehat{\psi}(x,y)(f(x) - f(y))^2 \\
&\leq 2 \sum_{x,y \in \Omega_i} \widehat{\psi}(x,y)[(f(x) - \mathrm{E}_{\pi_i} f)^2 + (f(y) - \mathrm{E}_{\pi_i} f)^2] \\
&\leq 2\eta \sum_{x \in A} \pi_i(x)(f(x) - \mathrm{E}_{\pi_i} f)^2 + 2\eta \sum_{y \in B} \pi_i(y)(f(y) - \mathrm{E}_{\pi_i} f)^2 \\
&\leq 2\eta \operatorname{Var}_{\pi_i} f,
\end{aligned}
$$

as required. $\square$

**4. Examples.** In this final section we collect a number of illustrative examples. The aim is more to rederive a variety of existing results in a simple, uniform manner than to obtain new results.

4.1. *Toy example.* Consider the symmetric random walk on the $2n$ vertex "pince-nez" graph in Figure 1 obtained by joining two disjoint $n$ cycles by a single edge. Suppose transitions within cycles occur with probability $\frac{1}{3}$, while the unique transition between cycles occurs with probability $p \leq \frac{1}{3}$. Loop probabilities are defined by complementation. The transition probabilities are symmetric, so the random walk is time-reversible and its stationary distribution is uniform.

Now decompose the set of vertices (states) into two disjoint subsets, $\Omega_0$ and $\Omega_1$, where $\Omega_0$ contains the $n$ vertices in the first cycle and $\Omega_1$ contains the $n$ vertices in the second cycle. The spectral gap for each cycle considered in isolation is $\frac{2}{3}(1-\cos(2\pi/n))$. (Diaconis and Saloff-Coste treated this example in [5], Section 4.2. The factor $\frac{2}{3}$ arises because our transition probabilities are $\frac{1}{3}$ instead of $\frac{1}{2}$.) Since $1-\cos x \geq 2x^2/5$ for $0 \leq x \leq \pi/2$, we have that the spectral gap for each restriction chain is at least $16\pi^2/15n^2$ (assuming $n \geq 4$), so we may take $\lambda_{\min} = 10n^{-2}$. The projection chain in this example is the symmetric two-state chain with transition probability $p/n$ between states, so we take $\bar{\lambda} = 2p/n$. Finally $\gamma = p$. (Recall that $\gamma$ is the maximum, over all states, of the probability of exiting from the current block of the partition of the states.) Theorem 1 yields, as Poincaré constant for the random walk on the pince-nez,

$$\lambda = \min\left\{\frac{2p}{3n}, \frac{20}{3n^3+2n^2}\right\}.$$

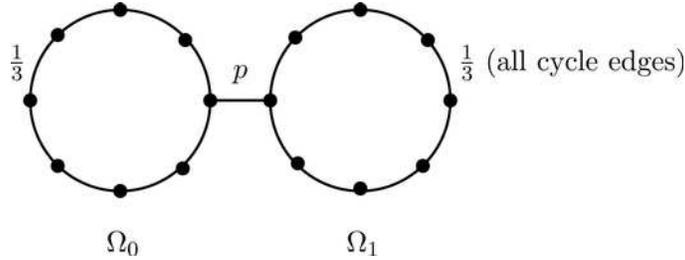

FIG. 1. *Pince-nez graph* $(n = 8)$.

Note that $\lambda = \Omega(n^{-3})$ when $p = \Omega(n^{-2})$ and $\lambda = \Omega(pn^{-1})$ when $p = O(n^{-2})$; in the latter case, our estimate is tight to within a constant factor, and a factor $n^2$ better than existing decomposition bounds which have the form $\lambda = \Theta(\bar{\lambda}\lambda_{\min}) = \Theta(pn^{-3})$.

A log-Sobolev inequality may be obtained equally simply by using a similar calculation. The restriction chains satisfy a log-Sobolev inequality with constant $16\pi^2/75n^2$ ([5], Section 4.2), so we may take $\alpha_{\min} = 2n^{-2}$. The log-Sobolev constant for the two-state projection chain is $\bar{\alpha} = p/n$ (see [5], Theorem A.2). Thus, by Theorem 4,

$$\alpha = \min\left\{\frac{p}{3n}, \frac{2}{3n^3+n^2}\right\}.$$

4.2. *Product of two Markov chains.* Consider two finite-state time-reversible Markov chains $(X, P_X)$ and $(Y, P_Y)$ with state spaces $X$ and $Y$. There are a number of ways to define a product Markov chain $(\Omega, P)$ on $\Omega = X \times Y$, but one which suits our purpose is to define the transition probabilities $P$ in terms of the transition probabilities $P_X$ and $P_Y$ as follows. For all $(x,y), (x',y') \in \Omega$,

$$P((x,y),(x',y')) := \begin{cases} P_X(x,x'), & \text{if } x' \neq x \text{ and } y' = y, \\ P_Y(y,y'), & \text{if } x' = x \text{ and } y' \neq y, \\ P_X(x,x) + P_Y(y,y) - 1, & \text{if } x' = x \text{ and } y' = y, \\ 0, & \text{otherwise.} \end{cases}$$

For the loop probabilities to be nonnegative, we require $P_X(x,x) + P_Y(y,y) \geq 1$ for all $(x,y) \in \Omega$ and we assume this from now on.

Our goal is to establish a Poincaré constant $\lambda$ for $(\Omega, P)$ in terms of those for $(X, P_X)$ and $(Y, P_Y)$: let us call them $\lambda_X$ and $\lambda_Y$, respectively.

For convenience, identify $X$ with $[n]$. Then, writing $\Omega_i := \{i\} \times Y$, we have the natural decomposition $\Omega = X \times Y = \bigcup_{i \in [n]} \Omega_i$. (Of course, we could have reversed the roles of $X$ and $Y$, and indexed the restriction chains by $Y$.) Each of the restriction chains is isomorphic to $(Y, P_Y)$ and so $\lambda_{\min} = \lambda_Y$. The projection chain is isomorphic to $(X, P_X)$ and so $\bar{\lambda} = \lambda_X$. By symmetry, $\hat{\pi}_i^j = \pi_i$ whenever the former is defined, and hence $\hat{\gamma} = 0$ and we are in the situation of Corollary 3. We obtain $\lambda = \min\{\lambda_X, \lambda_Y\}$ as the required Poincaré constant for $(\Omega, P)$, and this is tight. Exactly the same argument applies to the log-Sobolev constant.

4.3. *One-dimensional Ising model.* Consider the path of length $n$, that is, the graph with vertex set $[n]$ and edges joining vertices differing by 1. *Configurations* of the Ising model are just assignments $\sigma : [n] \to \{-1, +1\}$ of $\pm 1$ "spins" to the vertices of the graph. The *Hamiltonian* of the Ising system on the path is defined by

$$H(\sigma) := \sum_{i=0}^{n-2} [1 - \sigma(i)\sigma(i+1)]/2;$$

in other words, we count 1 for every pair of adjacent unlike spins. Denote the set of all $2^n$ configurations by $\Omega$. We wish to sample configurations from the Boltzmann–Gibbs distribution $\pi(\sigma) := \exp(-\beta H(\sigma))/Z$ on $\Omega$, where $Z := \sum_{\sigma \in \Omega} \exp(-\beta H(\sigma))$ is the *partition function* of the system and $\beta \in \mathbb{R}^+$ is *inverse temperature.* (What has been described is the *ferromagnetic* Ising model, which favors like spins; see [3] for background.)

One standard way to construct a Markov chain on $\Omega$ with stationary distribution $\pi$ is through single-site heat-bath dynamics. For $i \in [n]$ and $\sigma : [n] \to \{-1, +1\}$, let $\sigma_{[i \leftarrow +1]}$ (resp. $\sigma_{[i \leftarrow -1]}$) be the configuration that agrees with $\sigma$ at all vertices except possibly vertex $i$, where $\sigma_{[i \leftarrow +1]}(i) = +1$ [resp. $\sigma_{[i \leftarrow -1]}(i) = -1$]. The transitions of our heat-bath Markov chain are defined by the following trial, where $\sigma$ is the current state:

1. Select $i \in [n]$ uniformly at random.
2. Let

$$p := \frac{\exp\{-\beta H(\sigma_{[i \leftarrow +1]})\}}{\exp\{-\beta H(\sigma_{[i \leftarrow +1]})\} + \exp\{-\beta H(\sigma_{[i \leftarrow -1]})\}}.$$

Then with probability $p$, the new state is $\sigma_{[i \leftarrow +1]}$, and with probability $1 - p$, the new state is $\sigma_{[i \leftarrow -1]}$.

For convenience, we imagine that there are extra vertices 0 and $n$ with specified fixed spins, so that $p \in \{\frac{1}{2}, e^{\pm\beta}/(e^\beta + e^{-\beta})\}$.

Choose a vertex $m \in [n]$ as close to the midpoint of the path as possible (e.g., $m = \lfloor n/2 \rfloor$) and partition the configurations into two sets $\Omega = \Omega_+ \cup \Omega_-$,

BOUNDS ON POINCARÉ AND LOG-SOBOLEV CONSTANTS        17

where $\Omega_+$ (resp. $\Omega_-$) is the set of all configurations $\sigma$ with $\sigma(m) = +1$ [resp. $\sigma(m) = -1$]. Consider the restrictions of the Markov chain to $\Omega_+$ and $\Omega_-$, and the corresponding projection chain (which in this case has just two states).

A little optimization gives the spectral gap of the projection chain as $\bar{\lambda} \geq 1/(\cosh\beta)^2 n$. The parameter $\gamma$ satisfies $\gamma \leq 1/(1+e^{-2\beta})n$. Thus

$$\lambda \geq \min\left\{\frac{1}{3(\cosh\beta)^2 n}, \frac{\lambda_{\min}}{1+(3/4)(e^{2\beta}+1)}\right\}.$$

Each restriction chain is a direct product of two independent Ising systems on at most $\lfloor n/2 \rfloor$ vertices: independent because we fixed a spin at the middle of the path. The spectral gap of a direct product is the minimum of the spectral gaps of the components, as we saw in Section 4.2. So, denoting by $\lambda_k$ the spectral gap of the ferromagnetic Ising system on $[k]$ (with updates at any given site occurring at rate/probability $1/n$), we have the recurrence

$$\lambda_k \geq \min\left\{\frac{1}{3(\cosh\beta)^2 n}, \frac{\lambda_{\lfloor k/2 \rfloor}}{1+(3/4)(e^{2\beta}+1)}\right\}.$$

This has solution $\lambda_n = \Omega(n^{-c})$, where $c = 1 + \log_2\{1 + \frac{3}{4}(e^{2\beta}+1)\}$. So, at any fixed temperature, the spectral gap is bounded by an inverse polynomial in $n$ whose exponent $c$ tends to $1 + \log_2 \frac{5}{2} < 2.33$ as $\beta \to 0$. In light of (2), the mixing time scales as $n^{c+1}$.

A similar argument applies to the log-Sobolev constant. However, the bound on $\bar{\alpha}$, the log-Sobolev constant of the two-state projection chain, is $1/2(\cosh\beta)^2 n$ (see [5], Theorem A.2), which is worse than the bound we had for $\bar{\lambda}$ by a factor of 2. As a result, we obtain $\alpha_n = \Omega(n^{-c'})$, where $c' = 1 + \log_2\{1 + \frac{3}{2}(e^{2\beta}+1)\}$. So at any fixed temperature the log-Sobolev constant is bounded by an inverse polynomial in $n$ whose exponent $c'$ tends to $1 + \log_2 4 = 3$ as $\beta \to 0$. Although $c' > c$, we recall that, by (24), the mixing time scales as $n^{c'} \log n$, which is an improvement on $n^{c+1}$ for small enough $\beta$. In a sense, these are both poor results, since the one-dimensional Ising model does not exhibit a phase transition and we should expect mixing time $O(n \log n)$ at any temperature (although with a constant of proportionality depending on $\beta$). Note, however, that if we had used existing decomposition theorems, we would have lost a factor $n$ at each level of recursion, leading to a bound on spectral gap that diminishes with $n$ faster than any inverse polynomial.

When the temperature is sufficiently high (i.e., $\beta$ is sufficiently close to 0), we get a better bound by switching to Corollaries 2 or 5. Indeed, since $\hat{\gamma} \to 0$ as $\beta \to 0$, the bounds on both spectral gap and log-Sobolev constant are of the form $\Omega(1/n^{1+\delta})$ with $\delta$ tending to 0 as $n \to \infty$. (Very recent work [16] improved on this result by showing that $n^{1+\delta}$ may be replaced by $n$ (in the case of spectral gap) and $n \log n$ (in the case of log-Sobolev constant) for all $\beta$ up to some critical value. This remark applies equally to Section 4.4.)

4.4. *Ising and other spin models on trees.* The calculation of Section 4.3 carries over, with very little change to balanced trees of bounded degree. Thus we can treat balls of given radius in the so-called Bethe lattice of coordination number $r$; loosely, the infinite regular tree of degree $r$. Again, for fixed $r$ and $\beta$, the Poincaré and log-Sobolev constants are inverse polynomial. This was already known (see [10] for inverse polynomial bounds on spectral gap), although the fact that spin systems in general and the Ising model in particular have polynomial mixing time on trees is not so obvious.

4.5. *Walks on the Boolean cube and balanced matroids.* Consider the random walk on the $n$-dimensional Boolean cube in which transitions (including loops) all occur with probability $1/(n+1)$. There is a natural decomposition of the $n$-dimensional cube into two $(n-1)$-dimensional cubes connected by a perfect matching. Clearly $\hat{\gamma} = 0$ and we are in the situation of Corollary 3. Then $\bar{\lambda} = 2/(n+1)$ and, by induction, $\lambda = 2/(n+1)$. Likewise, $\bar{\alpha} = 1/(n+1)$ ([5], Theorem A.2, again) and hence $\alpha = 1/(n+1)$. These results are of course well known; see, for example, [5], Section 4.1.

More interestingly, there is a nontrivial Markov chain—the so-called bases-exchange walk on a balanced matroid—which retains just enough of the properties of the cube to be treatable by essentially the same approach. We shall be very brief in describing the setting; refer to [7] and [9] for a more expansive treatment. Let $E$ be a finite ground set and let $\mathcal{B} \subseteq 2^E$ a nonempty collection of subsets of $E$. We say that $\mathcal{B}$ forms the collection of *bases* of a *matroid* $M = (E, \mathcal{B})$ if the following "exchange axiom" holds: For every pair of bases $X, Y \in \mathcal{B}$ and every element $e \in Y \setminus X$, there exists an element $f \in X \setminus Y$ such that $X \cup \{e\} \setminus \{f\} \in \mathcal{B}$. It is easy to show that every basis must have the same size, which is called the *rank* of $M$. Denote by $m = |E|$ the size of the ground set and by $r$ the rank. A concrete example of a matroid is provided by the set of spanning trees of a finite undirected graph $G$. Here we interpret $E$ as the set of edges of $G$ and interpret $\mathcal{B}$ as the set of all spanning trees. The rank is of course $r = n - 1$, where $n$ is the number of vertices of $G$. The exchange axiom is easily checked.

Two absolutely central operations on matroids are contraction and deletion. If $e \in E$ is an element of the ground set of $M$, then the matroid $M \setminus e$ obtained by *deleting* $e$ has ground set $E \setminus \{e\}$ and bases $\mathcal{B}(M \setminus e) = \{X \subseteq E \setminus \{e\} : X \in \mathcal{B}(M)\}$; the matroid $M/e$ obtained by *contracting* $e$ has the same ground set as $M \setminus e$, but bases $\mathcal{B}(M/e) = \{X \subseteq E \setminus \{e\} : X \cup \{e\} \in \mathcal{B}(M)\}$.

The exchange axiom suggests a very natural Markov chain on the set of bases $\mathcal{B}$. Suppose the state (basis) at time $t$ is $X_t \in \mathcal{B}$. Then the state at time $t+1$ is obtained as the result of the following trial:

STEP 1. Choose $e \in E$ and $f \in X_t$ independently and uniformly at random.

STEP 2. If $X_t \cup \{e\} \setminus \{f\} \in \mathcal{B}$, then $X_{t+1} \leftarrow X_t \cup \{e\} \setminus \{f\}$ else $X_{t+1} \leftarrow X_t$.

It is straightforward to check, using the exchange axiom, that this Markov chain $(X_t)$ is irreducible; furthermore, it is clearly time-reversible and has uniform stationary distribution. It has been conjectured that $(X_t)$ is rapidly mixing (i.e., achieves total variation distance $\varepsilon > 0$ from the stationary distribution in a number of steps polynomial in $m$, $r$ and $\log \varepsilon^{-1}$) for all matroids $M$. Although there is little evidence in favor of this conjecture, there is at least an interesting class of matroids, namely the "balanced" matroids, for which rapid mixing was established [7]. The class of balanced matroids includes all regular matroids, and hence graphic matroids (i.e., ones whose set of bases may be realized as the set of all spanning trees of a graph).

Fix any $e \in E$. Since there is a natural isomorphism between $M$ and the disjoint union of $M \setminus e$ and $M/e$, the bases-exchange walk is a natural candidate for the decomposition method. Let $\Omega_0$ (resp. $\Omega_1$) be the set of bases that does not contain $e$ (resp. does contain $e$). Note that $\Omega_0$ is isomorphic to $M \setminus e$ and $\Omega_1$ is isomorphic to $M/e$, enabling us to argue inductively about the two restriction chains.

Rather than define the notion of balanced matroid explicitly here, let us just note two of the consequences of balance:

1. Contractions and deletions of a balanced matroid are themselves balanced, in particular, $M \setminus e$ and $M/e$ are balanced matroids.
2. The transitions of the bases-exchange walk that cross from $M \setminus e$ to $M/e$ and vice versa support a fractional matching. That is, there is a function $w : \Omega_0 \times \Omega_1 \to \mathbb{R}$ that satisfies:

- $w(x, y) \geq 0$ for all $x, y$;
- $w(x, y) > 0$ entails $P(x, y) > 0$;
- $\sum_y w(x, y) = \pi(\Omega_1)$ for all $x$ and $\sum_x w(x, y) = \pi(\Omega_0)$ for all $y$.

Define a new Markov chain on $\Omega$, with transition probabilities $\widehat{P}$, as follows. Transition probabilities within $\Omega_0$ and $\Omega_1$ are unchanged, so the restriction chains are also unchanged. Transitions from $x \in \Omega_0$ to $y \in \Omega_1$ and vice versa occur with probability $\widehat{P}(x, y) = \widehat{P}(y, x) := w(x, y)/(rm)$. Note that nonzero transition probabilities in the bases-exchange walk are all $1/(rm)$ and that the new Markov chain does not add any transitions to those already available. From these observations it follows immediately that for any pair of distinct states $x, y$, it is the case that $\widehat{P}(x, y) \leq P(x, y)$. Thus, it is enough to bound the Poincaré constant $\lambda$ for the Markov chain with modified transition probabilities $\widehat{P}$. (This is a trivial application of the comparison method in [6].)

The key point about $\widehat{P}$ is that it uses the fractional matching to spread the transitions between $\Omega_0$ and $\Omega_1$ evenly. It is easily checked that $\hat{\gamma} = 0$ and that we are in the situation of Corollary 3. The projection chain (derived from $\widehat{P}$ by projection onto $\{\Omega_0, \Omega_1\}$) has two states $\{0, 1\}$ with $\bar{\pi}(0) = \pi(\Omega_0)$ and $\bar{\pi}(1) = \pi(\Omega_1)$. Its transition probabilities are $\overline{P}(0,1) = \bar{\pi}(1)/(rm)$ and $\overline{P}(1,0) = \bar{\pi}(0)/(rm)$. A brief calculation establishes $\bar{\lambda} = 1/(rm)$. [This is directly from the definition (1) using $f(0) = \bar{\pi}(1)$ and $f(1) = -\bar{\pi}(0)$, a choice that is unique up to scaling among functions with expectation 0.] By induction, since $M \setminus e$ and $M/e$ are also balanced, $\lambda = 1/(rm)$.

Almost the same argument applies to the log-Sobolev constant $\alpha$. However, there is a small technicality that arises from the asymmetry of the projection chain. Specifically, we know only the two sides are balanced to the extent that

$$\frac{1}{m} \leq \frac{|\Omega_1|}{|\Omega_0|} \leq r.$$

Unfortunately, while the spectral gap is constant for the asymmetric two-state Markov chain, this is no longer the case for the log-Sobolev constant. More precisely, using [5], Theorem A.2, we have that

$$\bar{\alpha} = \Omega\left(\frac{1}{rm \ln(\pi(\Omega_0)^{-1} + \pi(\Omega_1)^{-1})}\right),$$

which worsens the bound by a factor $\ln(1/(r+m))$.

We can recover from this setback by noting that $w(x,y) \leq \min\{\pi(\Omega_0), \pi(\Omega_1)\}$ for all $x, y$, so that we could have defined $\widehat{P}$ by

$$\widehat{P}(x,y) = \widehat{P}(y,x) = \frac{w(x,y)}{rm \min\{\pi(\Omega_0), \pi(\Omega_1)\}}.$$

This definition boosts the transition probabilities between $\Omega_0$ and $\Omega_1$ when either $\pi(\Omega_0)$ or $\pi(\Omega_1)$ is small, while leaving all transition probabilities bounded by $1/(rm)$, as required. This modification to $\widehat{P}$ more than compensates for the effect just identified, and indeed the worst case (calculus) is when $\pi(\Omega_0) = \pi(\Omega_1) = \frac{1}{2}$. The bottom line is that we achieve $\alpha = 1/2(rm)$.

These bounds were recently obtained [9] using a longer and more complicated argument tailored to the specific application.

4.6. *Hard-core model on trees.* We conclude with an example where off-the-shelf decomposition results do not suffice. In this case, a reasonable bound on spectral gap may be deduced from Theorem 1. However, the problem that threatened to arise at the end of the previous example—namely that the log-Sobolev constant of a highly asymmetric two-state Markov chain may be arbitrarily close to 0—causes real problems here. Nevertheless, by tailoring our decomposition to the problem at hand, we obtain a reasonable bound on the log-Sobolev constant. (Just before the final version of this article was prepared, Martinelli, Sinclair and Weitz [16] announced substantially better bounds than those obtained here. They were able to show optimal, i.e., $O(n \log n)$, mixing over a range of fugacities.)

Let $T_d$ denote the tree, rooted at $v$, of depth $d$ and branching factor $\Delta$. (Thus, the degree of any vertex that is not the root or a leaf is $\Delta + 1$. The degree of the root $v$ itself is $\Delta$.) Consider the hard-core lattice gas model defined on the set $\Omega$ of independent sets of $T_d$. For a given fugacity $\lambda > 0$ (note that in this section only, $\lambda$ does *not* refer to spectral gap), we are interested in the Boltzmann–Gibbs distribution $\pi$ defined on $\Omega$, where

$$\pi(\sigma) \propto \lambda^{|\sigma|}$$

and $|\sigma|$ denotes the cardinality of the independent set $\sigma$.

Once again, a simple Markov chain with state space $\Omega$ and stationary distribution $\pi$ is the following single-site heat-bath dynamics known as the Glauber dynamics. For technical reasons we define the chain with respect to a parameter $N \geq n$, where $n$ denotes the number of vertices in $T_d$.

From $X_t \in \Omega$:

- Choose a vertex $z$ uniformly at random.
- Set

$$X' = \begin{cases} X_t \setminus \{z\}, & \text{with probability } 1/(1+\lambda), \\ X_t \cup \{z\}, & \text{with probability } \lambda/(1+\lambda). \end{cases}$$

- If $X' \in \Omega$, set $X_{t+1} = X'$ with probability $n/N$; otherwise, set $X_{t+1} = X_t$.

In practice, we set $N$ equal to the number $n$ of vertices in the tree at the top level of the inductive argument; the extra parameter eliminates rescalings in the upcoming analysis.

We decompose the state space as $\Omega = \Omega_1 \cup \Omega_2 \cup \Omega_3$, where

$$\Omega_1 := \{\sigma \in \Omega : v \in \sigma\},$$
$$\Omega_2 := \{\sigma \in \Omega \setminus \Omega_1 : \sigma \cup \{v\} \in \Omega\},$$
$$\Omega_3 := \Omega \setminus \Omega_1 \setminus \Omega_2.$$

Also let $\Omega_4 := \Omega_2 \cup \Omega_3$.

Without loss of generality, assume $\mathrm{E}_\pi f^2 = 1$. Using (25) twice,

$$\mathcal{L}_\pi(f) = \sum_{i=1}^{3} \bar{\pi}(i) \mathcal{L}_{\pi_i}(f) + \sum_{i=1}^{3} \bar{\pi}(i)(\mathrm{E}_{\pi_i} f^2) \ln(\mathrm{E}_{\pi_i} f^2)$$

(39)

$$\leq \bar{\pi}(1) \mathcal{L}_{\pi_1}(f) + \bar{\pi}(4) \mathcal{L}_{\pi_4}(f) + \sum_{i=1}^{3} \bar{\pi}(i)(\mathrm{E}_{\pi_i} f^2) \ln(\mathrm{E}_{\pi_i} f^2).$$

Observe that the restriction chain on $\Omega_1$ with transition probabilities $P_1$ (where $N$ is constant so that the transition probabilities are identical in $P_1$ and $P$) is simply the product chain of $\Delta^2$ copies of the chain on $T_{d-2}$. Similarly, the restriction chain on $\Omega_4$ is the product of $\Delta$ copies of the chain on $T_{d-1}$. The projection chain also has a simple structure. In particular, there is a bijection $\nu:\Omega_1 \to \Omega_2$ such that for $y = \nu(x)$ we have $P(x,y) = 1/(1+\lambda)N$ and $P(y,x) = \lambda/(1+\lambda)N$. This perfect matching $\nu$ captures the only transitions between the sets $\Omega_1$ and $\Omega_4$.

Define a new chain on state space $\{1,2,3\}$ with stationary distribution $\bar{\pi}$ (i.e., the same as that of the projection chain on $\{\Omega_1, \Omega_2, \Omega_3\}$) with transition probabilities $\widehat{P}$ given by

$$\widehat{P}(1,2) = \lambda/(1+\lambda),$$
$$\widehat{P}(2,1) = 1/(1+\lambda),$$
$$\widehat{P}(2,3) = \min\{1, \bar{\pi}(3)/\bar{\pi}(2)\},$$
$$\widehat{P}(3,2) = \min\{1, \bar{\pi}(2)/\bar{\pi}(3)\},$$
$$\widehat{P}(1,3) = \widehat{P}(3,1) = 0.$$

(The fictional Markov chain $\widehat{P}$ is a formal device: We establish and apply a log-Sobolev inequality for $\widehat{P}$ and then relate the various resulting terms to the actual chain $\overline{P}$.) Let $\hat{\alpha}$ denote the log-Sobolev constant of this chain. From (39) and the log-Sobolev inequality for $\widehat{P}$, we now have

(40)
$$\mathcal{L}_\pi(f) \leq \frac{\bar{\pi}(1)}{\alpha_{d-2}} \mathcal{E}_{\pi_1}(f,f) + \frac{\bar{\pi}(4)}{\alpha_{d-1}} \mathcal{E}_{\pi_4}(f,f)$$
$$+ \frac{1}{\hat{\alpha}} \sum_{i=1,2} \bar{\pi}(i)\widehat{P}(i, i+1)(\sqrt{\mathrm{E}_{\pi_i} f^2} - \sqrt{\mathrm{E}_{\pi_{i+1}} f^2})^2.$$

We need to bound the last terms on the right-hand side. Beginning with $i = 1$, observe that $\widehat{P}(1,2) = N\overline{P}(1,2)$. Moreover, the bijection $\nu$ implies the corresponding $\eta = 0$, and the following statement holds:

$$\bar{\pi}(1)\widehat{P}(1,2)(\sqrt{\mathrm{E}_{\pi_1} f^2} - \sqrt{\mathrm{E}_{\pi_2} f^2})^2 \leq N\mathcal{C}_{12}.$$

For $i = 2$ our bound proceeds as follows, starting with an inequality akin to (31):

$$\bar{\pi}(2)\widehat{P}(2,3)(\sqrt{\mathrm{E}_{\pi_2} f^2} - \sqrt{\mathrm{E}_{\pi_3} f^2})^2$$
$$\leq \bar{\pi}(2)\widehat{P}(2,3) \sum_{\substack{x \in \Omega_2 \\ y \in \Omega_3}} \pi_2(x)\pi_3(y)(f(x) - f(y))^2$$
$$\leq \min\{\bar{\pi}(2), \bar{\pi}(3)\} \operatorname{Var}_{\pi_4}(f)$$

$$\leq \frac{\min\{\bar{\pi}(2), \bar{\pi}(3)\}}{\alpha_{d-1}} \mathcal{E}_{\pi_4}(f, f).$$

Substituting these estimates into (40), we obtain

$$\mathcal{L}_\pi(f) \leq \frac{\bar{\pi}(1)}{\alpha_{d-2}} \mathcal{E}_{\pi_1}(f,f) + \frac{\bar{\pi}(4)}{\alpha_{d-1}} \left(1 + \frac{\min\{\bar{\pi}(2), \bar{\pi}(3)\}}{\bar{\pi}(4)\,\hat{\alpha}}\right) \mathcal{E}_{\pi_4}(f, f) + \frac{N\mathcal{C}_{12}}{\hat{\alpha}}$$

$$\leq \max\left\{\frac{1}{\alpha_{d-2}}, \left(1 + \frac{\min\{\bar{\pi}(2), \bar{\pi}(3)\}}{\bar{\pi}(4)\,\hat{\alpha}}\right)\frac{1}{\alpha_{d-1}}, \frac{N}{\hat{\alpha}}\right\} \mathcal{E}_\pi(f,f),$$

leading to the recurrence

(41) $$\alpha_d \geq \min\left\{\alpha_{d-2}, \left(1 + \frac{\min\{\bar{\pi}(2), \bar{\pi}(3)\}}{\bar{\pi}(4)\,\hat{\alpha}}\right)^{-1} \alpha_{d-1}, \frac{\hat{\alpha}}{N}\right\}.$$

To bound $\hat{\alpha}$, we use a result of Diaconis and Saloff-Coste ([5], Theorem A.1) on the log-Sobolev constant of the chain on $K_3$ (the complete graph on three vertices) with transition probabilities $P(i,j) = \bar{\pi}(j)$ for all $1 \leq i, j \leq 3$, that is,

(42) $$\alpha_{K_3} \geq \frac{1 - 2\bar{\pi}_{\min}}{\ln(1/\bar{\pi}_{\min} - 1)} \geq \frac{1}{3\ln(1/\bar{\pi}_{\min})},$$

where $\bar{\pi}_{\min} = \min\{\bar{\pi}(1), \bar{\pi}(2), \bar{\pi}(3)\}$. Comparing our three-state chain (with transition probabilities $\widehat{P}$) with $K_3$, we obtain the following claim.

CLAIM 7. *The quantities $\hat{\alpha}$ and $\alpha_{K_3}$ defined above are related by the inequality $\hat{\alpha} \geq \alpha_{K_3}/2(1 + \lambda)$.*

PROOF. For succinctness, write $\bar{f}(i) := \sqrt{\mathrm{E}_{\pi_i} f^2}$. Suppose $0 < \xi < 1$. Then

$$(\bar{f}(1) - \bar{f}(3))^2 \leq [(\bar{f}(1) - \bar{f}(2)) + (\bar{f}(2) - \bar{f}(3))]^2$$

$$= \left[\sqrt{\xi}\frac{1}{\sqrt{\xi}}(\bar{f}(1) - \bar{f}(2)) + \sqrt{1-\xi}\frac{1}{\sqrt{1-\xi}}(\bar{f}(2) - \bar{f}(3))\right]^2$$

$$\leq \frac{1}{\xi}(\bar{f}(1) - \bar{f}(2))^2 + \frac{1}{1-\xi}(\bar{f}(2) - \bar{f}(3))^2,$$

where the inequality is Cauchy–Schwarz. Now set

$$\xi = \frac{(\bar{\pi}(1) + \bar{\pi}(2))\bar{\pi}(3)}{\bar{\pi}(1)\bar{\pi}(2) + 2\bar{\pi}(1)\bar{\pi}(3) + \bar{\pi}(2)\bar{\pi}(3)}.$$

Then

$$\bar{\pi}(1)\bar{\pi}(3)(\bar{f}(1) - \bar{f}(3))^2$$

$$\leq \left(\bar{\pi}(1) + \frac{2\bar{\pi}(1)\bar{\pi}(3)}{\bar{\pi}(2)} + \bar{\pi}(3)\right)$$

$$\times \left[\frac{\bar{\pi}(1)\bar{\pi}(2)}{\bar{\pi}(1) + \bar{\pi}(2)}(\bar{f}(1) - \bar{f}(2))^2 + \frac{\bar{\pi}(2)\bar{\pi}(3)}{\bar{\pi}(2) + \bar{\pi}(3)}(\bar{f}(2) - \bar{f}(3))^2\right].$$

Starting with the log-Sobolev inequality for the Diaconis–Saloff-Coste chain on $K_3$ and applying the above inequality, we obtain

$$\alpha_{K_3}\mathcal{L}_{\bar{\pi}}(\bar{f}) \leq \bar{\pi}(1)\bar{\pi}(2)(\bar{f}(1) - \bar{f}(2))^2$$
$$+ \bar{\pi}(1)\bar{\pi}(3)(\bar{f}(1) - \bar{f}(3))^2 + \bar{\pi}(2)\bar{\pi}(3)(\bar{f}(2) - \bar{f}(3))^2$$
$$\leq 2\left(1 + \frac{\bar{\pi}(1)\bar{\pi}(3)}{\bar{\pi}(2)}\right)\left[\frac{\bar{\pi}(1)\bar{\pi}(2)}{\bar{\pi}(1) + \bar{\pi}(2)}(\bar{f}(1) - \bar{f}(2))^2\right.$$
$$\left. + \frac{\bar{\pi}(2)\bar{\pi}(3)}{\bar{\pi}(2) + \bar{\pi}(3)}(\bar{f}(2) - \bar{f}(3))^2\right]$$
$$\leq 2\left(1 + \frac{\bar{\pi}(1)\bar{\pi}(3)}{\bar{\pi}(2)}\right)\mathcal{E}_{\bar{\pi}}(\bar{f}, \bar{f}),$$

where $\mathcal{E}_{\bar{\pi}}(\bar{f}, \bar{f})$ is the Dirichlet form associated with the $\widehat{P}$ chain. Hence

$$\hat{\alpha} \geq \frac{\alpha_{K_3}}{2(1 + \bar{\pi}(1)\bar{\pi}(3)/\bar{\pi}(2))} \geq \frac{\alpha_{K_3}}{2(1 + \lambda)},$$

as claimed, where the second inequality uses the fact that $\bar{\pi}(1) = \lambda\bar{\pi}(2)$. □

Combined with our earlier estimate (42) for $\alpha_{K_3}$, Claim 7 gives

$$\hat{\alpha} \geq \frac{1}{6(1+\lambda)\ln(1/\bar{\pi}_{\min})}.$$

By considering appropriate mappings $\Omega \to \Omega_i$, we see that

$$\bar{\pi}(1) \geq \frac{\lambda}{(1+\lambda)^{\Delta+1}}, \qquad \bar{\pi}(2) \geq \frac{1}{(1+\lambda)^{\Delta+1}}, \qquad \bar{\pi}(3) \geq \frac{\lambda}{(1+\lambda)^{\Delta+2}}.$$

[For example, for the first of these, consider the mapping that forces $\sigma \mapsto \sigma \cup \{v\} \setminus \Gamma(v)$, where $\Gamma(v)$ denotes the set of neighbors of $v$.] Thus

$$\frac{1}{\bar{\pi}_{\min}} \leq \frac{(1+\lambda)^{\Delta+2}}{\min\{1,\lambda\}}.$$

In the other direction, similar arguments yield

$$\bar{\pi}(2) \leq \frac{1}{\lambda}, \qquad \bar{\pi}(3) \leq (1+\lambda)^\Delta - 1, \qquad \bar{\pi}(4) \geq \frac{1}{1+\lambda},$$

so that

$$\min\{\bar{\pi}(2), \bar{\pi}(3)\} \leq \min\left\{\frac{1}{\lambda}, (1+\lambda)^\Delta - 1, 1\right\}.$$

Combining these various inequalities yields

$$\frac{\min\{\bar{\pi}(2), \bar{\pi}(3)\}}{\bar{\pi}(4)\hat{\alpha}} \leq \frac{6\min\{\bar{\pi}(2), \bar{\pi}(3)\}(1+\lambda)\ln(1/\bar{\pi}_{\min})}{\bar{\pi}(4)}$$
$$\leq g_\Delta(\lambda),$$

where

(43) $\quad g_\Delta(\lambda) = 6 \min\left\{\dfrac{1}{\lambda}, (1+\lambda)^\Delta - 1, 1\right\}(1+\lambda)^2 \ln\left(\dfrac{(1+\lambda)^{\Delta+2}}{\min\{1,\lambda\}}\right).$

Returning to recurrence (41),

$$1/\alpha_d \leq \max\{(1 + g_\Delta(\lambda))/\alpha_{d-1}, N/\hat{\alpha}\},$$

leading to

$$\alpha^{-1} \leq n(1 + g_\Delta(\lambda))^{\log_\Delta n}$$
$$= n^{1+\log_\Delta(1+g_\Delta(\lambda))}.$$

To make sense of this bound on the log-Sobolev constant $\alpha$, observe that

$$g_\Delta(\lambda) = O(\Delta\lambda(1 + |\log \lambda|)).$$

This estimate can be obtained from (43) by considering separately the ranges $\lambda < \Delta^{-1}$, $\Delta^{-1} \leq \lambda < 1$ and $\lambda \geq 1$. In particular, we see that $g_\Delta(\lambda) \to 0$ as $\lambda \to 0$ for any fixed $\Delta$. In rough terms, the mixing time tends to be linear as the fugacity $\lambda$ tends to 0.

## REFERENCES

[1] CARACCIOLO, S., PELISSETTO, A. and SOKAL, A. D. (1992). Two remarks on simulated tempering. Unpublished manuscript.
[2] CESI, F. (2001). Quasi-factorization of the entropy and logarithmic Sobolev inequalities for Gibbs random fields. *Probab. Theory Related Fields* **120** 569–584. MR1853483
[3] CIPRA, B. A. (1987). An introduction to the Ising model. *Amer. Math. Monthly* **94** 937–959. MR936054
[4] COOPER, C., DYER, M. E., FRIEZE, A. M. and RUE, R. (2000). Mixing properties of the Swendsen–Wang process on the complete graph and narrow grids. *J. Math. Phys.* **41** 1499–1527. MR1757967
[5] DIACONIS, P. and SALOFF-COSTE, L. (1996). Logarithmic Sobolev inequalities for finite Markov chains. *Ann. Appl. Probab.* **6** 695–750. MR1410112
[6] DIACONIS, P. and SALOFF-COSTE, L. (1993). Comparison theorems for reversible Markov chains. *Ann. Appl. Probab.* **3** 696–730. MR1233621
[7] FEDER, T. and MIHAIL, M. (1992). Balanced matroids. In *Proceedings of the 24th Annual ACM Symposium on Theory of Computing* 26–38. ACM Press, New York.
[8] JERRUM, M. (2003). *Counting, Sampling and Integrating: Algorithms and Complexity*. Birkhäuser, Basel. MR1960003
[9] JERRUM, M. and SON, J.-B. (2002). Spectral gap and log-Sobolev constant for balanced matroids. In *Proceedings of the 43rd IEEE Annual Symposium on Foundations of Computer Science* 721–729. Computer Society Press, New York.
[10] KENYON, C., MOSSEL, E. and PERES, Y. (2001). Glauber dynamics on trees and hyperbolic graphs. In *Proceedings of the 42nd IEEE Annual Symposium on Foundations of Computer Science* 568–578. Computer Society Press, New York. MR1948746

[11] KIPNIS, C., OLLA, S. and VARADHAN, S. R. S. (1989). Hydrodynamics and large deviation for simple exclusion processes. *Comm. Pure Appl. Math.* **42** 115–137. MR978701
[12] LEE, T.-Y. and YAU, H.-T. (1998). Logarithmic Sobolev inequality for some models of random walks. *Ann. Probab.* **26** 1855–1873. MR1675008
[13] MADRAS, N. and RANDALL, D. (2002). Markov chain decomposition for convergence rate analysis. *Ann. Appl. Probab.* **12** 581–606. MR1910641
[14] MARTIN, R. A. (2001). Paths, sampling, and Markov chain decomposition. Ph.D. dissertation, Georgia Institute of Technology.
[15] MARTIN, R. A. and RANDALL, D. (2000). Sampling adsorbing staircase walks using a new Markov chain decomposition method. In *Proceedings of the 41st IEEE Annual Symposium on Foundations of Computer Science* 492–502. Computer Society Press, New York. MR1931846
[16] MARTINELLI, F., SINCLAIR, A. and WEITZ, D. (2003). The Ising model on trees: Boundary conditions and mixing time. Technical Report CSD-03-1256, Computer Science Division, Univ. California, Berkeley. (Extended abstract available in *Proceedings of the 44th Annual IEEE Symposium on Foundations of Computer Science.* Computer Society Press, New York.)

M. JERRUM
J.-B. SON
SCHOOL OF INFORMATICS
UNIVERSITY OF EDINBURGH
THE KING'S BUILDINGS
EDINBURGH EH9 3JZ
UNITED KINGDOM
E-MAIL: mrj@inf.ed.ac.uk

P. TETALI
SCHOOL OF MATHEMATICS
GEORGIA INSTITUTE OF TECHNOLOGY
ATLANTA, GEORGIA 30332-0160
USA

E. VIGODA
DEPARTMENT OF COMPUTER SCIENCE
UNIVERSITY OF CHICAGO
1100 E. 58TH STREET
CHICAGO, ILLINOIS 60637
USA